\documentclass[11pt,reqno]{amsart}

\usepackage{iftex}
\usepackage{amsmath}
\ifPDFTeX
  \usepackage[T1]{fontenc}
  \usepackage{lmodern,amssymb}
\else
  \usepackage{fontspec,unicode-math}
\fi
\usepackage[a4paper,margin=30mm]{geometry}
\usepackage{microtype}
\usepackage{booktabs}
\usepackage[hidelinks]{hyperref}

\newtheorem{theorem}{Theorem}[section]
\newtheorem{proposition}[theorem]{Proposition}
\newtheorem{lemma}[theorem]{Lemma}
\newtheorem{corollary}[theorem]{Corollary}
\theoremstyle{definition}
\newtheorem{definition}[theorem]{Definition}

\newcommand{\R}{\mathbb R}
\newcommand{\T}{\mathbb T}
\newcommand{\E}{\mathbb E}
\newcommand{\F}{\mathcal F}
\newcommand{\Rank}{\mathcal R}
\newcommand{\Span}{\operatorname{span}}
\newcommand{\norm}[1]{\lVert #1\rVert}
\newcommand{\abs}[1]{\lvert #1\rvert}
\allowdisplaybreaks[1]

\title[Dunford--Pettis property of Sobolev martingales]%
{The Dunford--Pettis property and complemented Hilbert subspaces of Sobolev martingale spaces}
\author{Juliusz Pham}
\date{}
\hypersetup{pdftitle={The Dunford--Pettis property and complemented Hilbert subspaces of Sobolev martingale spaces},pdfauthor={Juliusz Pham}}

\begin{document}

\begin{abstract}
Let \(L^1_W\) be the space of integrable terminal values of martingales on a homogeneous tree whose local differences belong to a fixed matrix space \(W\). We show that \(L^1_W\) has the Dunford--Pettis property if and only if \(W\) is spanned by its rank-one matrices. Every complemented reflexive subspace of \(L^1_W\) is isomorphic to a Hilbert space. Under the second structural condition of Ayoush, Stolyarov and Wojciechowski, the martingale fractional integral of order \(1/2\) is completely continuous and fails to be absolutely \(2\)-summing. This gives nonisomorphism with \(L^1\) for \(W\ne0\). On the \(4\)-ary tree we exhibit all three possible combinations of the Dunford--Pettis property and the existence of a complemented copy of \(\ell^2\). In one example a bounded operator into \(\ell^4\) fails to be completely continuous, although every bounded operator into a Hilbert space is completely continuous. We also prove that, for every \(n\ge1\), each complemented reflexive subspace of \(W^{1,1}(\T^n)\) is isomorphic to a Hilbert space.
\end{abstract}

\maketitle

\section{Introduction}

On the torus, a vector field is a gradient precisely when it has mean
zero and its Fourier coefficient at each nonzero frequency \(\xi\) is
parallel to \(\xi\). Thus the gradient identifies the mean-zero part of
\(W^{1,1}(\T^n)\), equipped with the gradient norm, with a closed
subspace of \(L^1(\T^n;\R^n)\). The constraint on the field matters for
the Banach-space structure. In dimension two, for example, Kislyakov
used the critical Sobolev embedding to distinguish \(W^{1,1}(\T^2)\)
from \(L^1\) \cite{Kislyakov}.

Ayoush, Stolyarov and Wojciechowski \cite{ASW} introduced a martingale
analogue of such constrained spaces. At each vertex of a homogeneous
tree, the values of a martingale difference on the children form a
matrix. One fixes a linear subspace \(W\) of admissible matrices and
imposes the same restriction at every vertex. Their second structural
condition ensures an endpoint estimate for the martingale fractional
integral. We use this estimate to examine the isomorphism type of
\(L^1_W\), the space of integrable terminal values satisfying the
constraint. The role of the Dunford--Pettis property in these models
was also raised in \cite{ASWRank}.

Which local constraints on the increments preserve the Dunford--Pettis
property? Which permit a complemented copy of a Hilbert space?
These questions come apart in a two-coordinate example. One part of
the argument transfers to classical Sobolev spaces: the approximation
obtained by stopping martingales has a counterpart via Lipschitz
truncation. This gives a structure theorem for every complemented
reflexive subspace of \(W^{1,1}(\T^n)\) in every dimension
(Theorem~\ref{thm:torus}).

\subsection{The Dunford--Pettis property}

A bounded operator is \emph{completely continuous} if it maps weakly null sequences to norm-null sequences. A Banach space \(X\) has the \emph{Dunford--Pettis property} (DPP) if every weakly compact operator \(X\to Z\), for any Banach space \(Z\), is completely continuous. Classical positive examples are \(L^1\) by Dunford and Pettis \cite{DP} and \(C(K)\) for compact Hausdorff \(K\) by Grothendieck \cite{GrothC}.

Bourgain \cite{Bourgain} established the property for ball and polydisc algebras and for spaces \(C^k([0,1]^d)\) equipped with the uniform norm on their derivatives of order at most \(k\). For related results on rich subspaces of \(C(K)\) and uniform algebras, see \cite{Lefevre,Saccone}.

The DPP passes to complemented subspaces. Indeed, a weakly compact operator on a complemented subspace extends to a weakly compact operator on the whole space by composition with the projection. The identity on an infinite-dimensional reflexive space is weakly compact and cannot be completely continuous, by weak sequential compactness of its unit ball and the failure of norm compactness. Thus a space with the DPP has no infinite-dimensional complemented reflexive subspace.

\section{The spaces and the structural conditions}
\label{sec:model}

Fix integers \(m\ge 2\) and \(k\ge 1\). Let \(\Omega=\{1,\ldots,m\}^{\mathbb N}\) with its uniform product probability measure \(\mu\), and let \(\F_n\) be the sigma algebra of the first \(n\) coordinates, with \(\F_0\) trivial. Write \(\E_n\) for conditional expectation and \(d_nf=\E_nf-\E_{n-1}f\). For a parent cylinder \(Q\) of generation \(n-1\), let \(Q_1,\ldots,Q_m\) be its children in the order of the \(n\)-th digit. Define
\[
 \begin{split}
 V_m&=\left\{v\in\R^m:\ \sum_{i=1}^m v_i=0\right\},\\
 D_Q f&=\bigl(d_nf|_{Q_1},\ldots,d_nf|_{Q_m}\bigr)
          \in V_m\otimes\R^k .
 \end{split}
\]
We identify \(V_m\otimes\R^k\) with matrices whose rows in \(\R^k\) sum to zero. In particular, \(v\otimes a=(v_i a)_{i=1}^m\).

\begin{definition}
For a linear subspace \(W\subset V_m\otimes\R^k\) and \(p=1,2\), put
\[
 L^p_W=\{f\in L^p(\Omega;\R^k): D_Q f\in W\text{ for every parent }Q\},
\]
with the ambient \(L^p\) norm. The value at the root, \(\E f\), is unrestricted. The term \emph{block} on \((a,b]\) denotes a member \(u\) with \(\E_a u=0\) and \(\E_bu=u\).
\end{definition}

The defining constraints are closed. The maps \(\E_n\) are finite-rank contractions on \(L^p_W\), with \(\E_nf\to f\) in \(L^p\). Multiplying a block by a bounded \(\F_a\)-measurable scalar preserves its local differences in \(W\) and its block support in levels. This follows by taking the scalar outside conditional expectations on \((a,b]\).

Following \cite{ASW}, we use two structural conditions. Here \(e_i\) denotes the \(i\)-th coordinate vector of \(\R^m\) and \(\mathbf 1=(1,\ldots,1)\).
\begin{align*}
 \tag{S1}
 &W\text{ contains no nonzero }v\otimes a,\quad v\in V_m,\ a\in\R^k.\\
 \tag{S2}
 &W\text{ contains no nonzero }(me_i-\mathbf 1)\otimes a,\quad
                 1\le i\le m,\ a\in\R^k.
\end{align*}
The first condition excludes all rank-one directions. The second excludes the particular directions of a martingale concentrating on one child. The second follows from the first.

Under (S1), Janson's theorem \cite{Janson} (see also \cite[Theorem~1.8]{ASW}) gives the martingale Hardy space estimate
\begin{equation}
 \label{eq:Janson}
 \E\left(\sum_{n\ge1}\abs{d_nf}^{2}\right)^{1/2}
 \lesssim_W \norm{f}_1,\qquad f\in L^1_W .
\end{equation}
Together with the reverse square-function bound, this yields the following standard consequence of Janson's theorem, familiar as folklore: for \(W\ne0\),
\[
 L^1_W\simeq H^1_{\mathrm{dyad}}\simeq \operatorname{Re} H^1(\T)
\]
as real Banach spaces. Indeed, with the root value included in the Hardy norm, Janson's estimate identifies \(L^1_W\) with the Hardy-norm space with constraint \(W\). The isomorphism of this constrained Hardy space with scalar martingale \(H^1\) is recorded in \cite[Theorem~6.4]{KL}. Maurey's theorems identify scalar \(m\)-adic and dyadic \(H^1\), and identify the latter with real Hardy space on the circle \cite[Theorems~4.14 and~5.4]{Maurey}. Here \(\operatorname{Re} H^1(\T)\) carries the norm \(\norm{u}_1+\norm{\mathcal H u}_1\), with \(\mathcal H\) the circular Hilbert transform. For \(W=0\), \(L^1_W\) consists of constants.

Under (S2), the martingale Riesz potential (fractional integral)
\[
 I_\alpha f=\sum_{n\ge1}m^{-\alpha n}d_nf,\qquad \alpha=1-1/p,
\]
maps \(L^1_W\) boundedly into the Lorentz space \(L^{p,1}(\Omega;\R^k)\) for \(1<p<\infty\). Ayoush, Stolyarov and Wojciechowski \cite[Theorem~1.9 and (1.6)]{ASW} prove the stronger estimate
\begin{equation}
 \label{eq:ASW}
 \sum_{n\ge1}m^{-(1-1/p)n}\norm{d_nf}_{L^{p,1}}
 \lesssim_{W,p}\norm{f}_1 .
\end{equation}
Their framework treats bounded-in-\(L^1\) martingales as boundary measures and also gives trace estimates and lower bounds on the Hausdorff dimension of these measures \cite{ASW,StolyarovTrace}. We work with martingales of the form \((\E_nf)\), so that the norm is that of the integrable terminal value \(f\).

\subsection{Main results}

\begin{enumerate}
\item \textbf{The Dunford--Pettis property.} The space \(L^1_W\) has
the DPP precisely when \(W\) is spanned by its rank-one matrices
(Theorem~\ref{thm:criterion}).
\item \textbf{A Jordan-type constraint.} For a specific two-coordinate constraint, a bounded operator into \(\ell^4\) fails to be completely continuous, while every bounded operator into a Hilbert space is completely continuous (Theorem~\ref{thm:jordan}).
\item \textbf{Complemented reflexive subspaces.} Every complemented reflexive
subspace of \(L^1_W\) or \(W^{1,1}(\T^n)\), \(n\ge1\), is isomorphic to a Hilbert space (Theorems~\ref{thm:hilbert} and~\ref{thm:torus}). In \(L^1_W\), a complemented copy of \(\ell^2\) exists precisely when some bounded operator into \(\ell^2\) fails to be completely continuous.
\item \textbf{The fractional integral.} For nonzero \(W\) satisfying
(S2), the bounded operator \(I_{1/2}:L^1_W\to L^2\) fails to be
absolutely \(2\)-summing, so \(L^1_W\) is not isomorphic to \(L^1\)
(Proposition~\ref{prop:noniso}). Nonisomorphism can also occur when
(S2) fails.
\end{enumerate}

\section{The fractional integral and nonisomorphism with \texorpdfstring{\(L^1\)}{L1}}
\label{sec:summing}

An operator \(A:X\to Y\) is \emph{absolutely \(2\)-summing} if there is a constant \(C\) such that for all finite families \((x_j)\) in \(X\),
\begin{equation}
 \label{eq:2summing}
 \left(\sum_j\norm{Ax_j}^{2}\right)^{1/2}
 \le C\sup_{x^*\in B_{X^*}}
            \left(\sum_j \abs{x^*(x_j)}^2\right)^{1/2}.
\end{equation}
Grothendieck's theorem asserts that every bounded operator from an \(L^1\) space to a Hilbert space is absolutely \(1\)-summing, and hence \(2\)-summing \cite{DJT,KPW}. Kislyakov used the failure of \(2\)-summability of the critical Sobolev embedding to distinguish \(W^{1,1}(\T^2)\) from \(L^1\) \cite{Kislyakov}. For the subspaces of \(L^1\) generated by the even or odd levels of the Haar system, the nonisomorphism with \(L^1\) was proved by Speckhofer \cite{Speckhofer}. The martingale fractional integral supplies the corresponding operator for every nonzero \(W\) satisfying (S2).

\begin{proposition}
 \label{prop:noniso}
 If \(W\ne 0\) satisfies \textup{(S2)}, then \(I_{1/2}:L^1_W\to L^2(\Omega;\R^k)\) is bounded and fails to be absolutely \(2\)-summing. In particular, \(L^1_W\) is not isomorphic to an \(L^1\) space.
\end{proposition}

\begin{proof}
 Boundedness follows from \eqref{eq:ASW} and the inclusion \(L^{2,1}\hookrightarrow L^2\). Choose \(w=(w_i)_{i=1}^m\in W\) with \(m^{-1}\sum_i\abs{w_i}^2=1\). For each cylinder \(Q\) of generation \(n-1\), define \(h_Q\) by \(\mu(Q)^{-1/2} w_i\) on \(Q_i\) and by zero off \(Q\). The functions \(h_Q\) are in \(L^2_W\) and are orthonormal in \(L^2\), by disjoint supports at each generation and orthogonality of martingale differences across generations.

 For each \(x^*\in B_{(L^1_W)^*}\), its restriction to \(L^2_W\) has norm at most one since \(\norm{f}_1\le\norm{f}_2\). Bessel's inequality shows that the supremum in \eqref{eq:2summing} is at most one for every finite family of the \(h_Q\). Each \(h_Q\) of generation \(n-1\) has \(I_{1/2}h_Q=m^{-n/2}h_Q\). Summing over all parents through generation \(N-1\) gives
 \[
   \sum_{n=1}^{N}\sum_{Q\in\F_{n-1}\,\text{atoms}}\norm{I_{1/2}h_Q}_2^2
       =\sum_{n=1}^{N}m^{n-1}m^{-n}=N/m .
 \]
 This violates \eqref{eq:2summing}. The isomorphism assertion follows because the property that all bounded operators into a Hilbert space are \(2\)-summing is preserved by isomorphism of the domain.
\end{proof}

The converse of the nonisomorphism assertion fails. The triadic example of Kazaniecki and Lechner \cite[Theorem~6.6]{KL} has a scalar summand with all mean-zero child patterns allowed, so (S2) fails. Its other summand is isomorphic to \(H^1\) and contains a complemented copy of \(\ell^2\), so the space fails the DPP and cannot be isomorphic to \(L^1\).
The scalar constraint \(W=V_4\) also fails (S2) and gives
all of \(L^1(\Omega)\). Thus failure of (S2) alone does not determine whether \(L^1_W\) is isomorphic to \(L^1\).

For the boundedness of \(I_{1/2}\) itself, (S2) is also necessary.
Indeed, suppose \((me_i-\mathbf1)\otimes a\in W\) with \(|a|=1\),
and let \(Q_N\) be the cylinder whose first \(N\) digits are all \(i\).
The function \(f_N=m^N\mathbf1_{Q_N}a\) has \(L^1\)-norm one and
belongs to \(L^1_W\). For \(1\le n\le N\), its only nonzero local
difference at level \(n\) is on the parent \(Q_{n-1}\), where its
pattern is \(m^{n-1}(me_i-\mathbf1)\otimes a\). In particular,
\[
    \|d_nf_N\|_2^2=(m-1)m^{n-1},
    \qquad
    \|I_{1/2}f_N\|_2^2=\frac{m-1}{m}N
\]
by orthogonality. The operator is consequently unbounded on these
finite martingales.

\section{Weakly compact sets and complemented reflexive subspaces}
\label{sec:complement}

\begin{lemma}[Martingale stopping]
 \label{lem:stop}
 Let \(K\subset L^1_W\) be bounded and uniformly integrable. For every \(\epsilon>0\) there is \(M<\infty\) such that each \(f\in K\) is within \(\epsilon\) in \(L^1\) of a \(g\in L^1_W\) with \(\norm{g}_\infty\le M\). If \(f\) is a block, \(g\) can be chosen on the same interval of levels.
\end{lemma}

\begin{proof}
 Take \(R\) larger than \(\sup_{f\in K}\norm f_1\). Stop \((\E_n f)\) when \(\E_n\abs f>R\) for the first time. In the block case start monitoring at the first endpoint \(a\) and ignore times after the last endpoint \(b\). In that case a stop at time \(a\) has stopped value \(\E_a f=0\). At a later crossing, \(\E_n\abs f\le mR\), since the value of this nonnegative martingale on a child is at most \(m\) times that on its parent. The terminal value \(g\) of the stopped martingale is therefore bounded by \(mR\). Each difference is the original difference multiplied by the indicator that it has yet to be stopped at its parent. Thus \(g\) has the required constraint, and in the block case \(\E_ag=0\) and \(\E_bg=g\).

 Let \(D\) be the union of the stopping cylinders. The maximal inequality gives \(\mu(D)\le \norm f_1/R\). Off \(D\), \(g=f\), and on \(D\) the inequality \(\int_D\abs g\le\int_D\abs f\) follows from the formula for the stopped value on each cylinder. The error is bounded by \(2\int_D\abs f\), which tends to zero uniformly in \(f\in K\) as \(R\to\infty\).
\end{proof}

\begin{lemma}[Approximation through a Hilbert space]
\label{lem:hilbert-approx}
Let \(X\) be a Banach space and \(i:H\to X\) a bounded map from a Hilbert space. Suppose that for every weakly compact \(K\subset X\) and every \(\epsilon>0\), some \(M<\infty\) satisfies
\begin{equation}
 \label{eq:approx}
 K\subset M i(B_H)+\epsilon B_X.
\end{equation}
Then a bounded operator \(A:X\to Z\) into a Banach space is completely continuous if and only if \(Ai\) is compact. Every complemented reflexive subspace of \(X\) is isomorphic to a Hilbert space.
\end{lemma}

\begin{proof}
If \(Ai\) is compact, \eqref{eq:approx} makes \(A(K)\) totally bounded. Applying this to a weakly convergent sequence together with its limit gives complete continuity. Conversely, \(i(B_H)\) is weakly compact, and a completely continuous operator maps it to a relatively norm-compact set by the Eberlein--\v{S}mulian theorem.

Let \(P:X\to Y\) be a bounded projection onto a reflexive subspace. Apply \eqref{eq:approx} to \(B_Y\) and take \(\epsilon\norm P\le 1/2\). Applying \(P\) gives \(B_Y\subset M Pi(B_H)+\tfrac12 B_Y\). Iterating the approximation of the residual expresses each \(y\in B_Y\) as \(Pih\) with \(\norm h\le 2M\). Thus \(Pi\) is surjective and \(Y\) is isomorphic to the Hilbert-space quotient \(H/\ker(Pi)\).
\end{proof}

\begin{theorem}
 \label{thm:hilbert}
 For every fixed \(W\), each complemented reflexive subspace of \(L^1_W\) is isomorphic to a Hilbert space. Moreover, with \(i:L^2_W\hookrightarrow L^1_W\) the natural inclusion, the following assertions hold.
 \begin{enumerate}
 \item For every Banach space \(Z\) and bounded \(A:L^1_W\to Z\), \(A\) is completely continuous if and only if \(Ai\) is compact.
 \item \(L^1_W\) contains a complemented copy of \(\ell^2\) if and only if there exists a bounded operator \(A:L^1_W\to\ell^2\) that is not completely continuous.
 \end{enumerate}
\end{theorem}

\begin{proof}
 Write \(X=L^1_W\) and \(H=L^2_W\). By the Dunford--Pettis weak compactness criterion in \(L^1\), every weakly compact subset of \(X\) is uniformly integrable. Lemma~\ref{lem:stop} gives \eqref{eq:approx}, so Lemma~\ref{lem:hilbert-approx} gives the Hilbertian conclusion and the criterion for complete continuity.

 Suppose \(A:X\to\ell^2\) is not completely continuous. Then \(Ai:H\to\ell^2\) is noncompact. By the spectral theorem applied to \((Ai)^*Ai\), there is an infinite-dimensional closed subspace \(H_0\) of \(H\) on which \(Ai\) is bounded below. On \(H_0\) the inclusion \(i\) is then bounded below as well, because \(\norm{Aih}_2\le\norm A\norm{ih}_1\). Let \(Y=i(H_0)\), and let \(Q\) be the orthogonal projection of \(\ell^2\) onto \(A(Y)\). The map \((A|_Y)^{-1}QA\) is a bounded projection onto \(Y\simeq\ell^2\). In the opposite direction, a projection onto a copy of \(\ell^2\), followed by an isomorphism onto \(\ell^2\), fails complete continuity on the image of its unit vector basis.

\end{proof}

The Hilbertian conclusion is an analogue of \cite[Corollary 2.1]{KP} for analytic \(H^1\). The proof here uses a uniform approximation on weakly compact sets, which is available for every \(W\) and does not require either structural condition.

\begin{corollary}
 \label{cor:Icc}
 If \(W\) satisfies \textup{(S2)}, then \(I_{1/2}:L^1_W\to L^2\) is completely continuous. If \(W\ne0\), this operator is noncompact.
\end{corollary}

\begin{proof}
 On \(L^2_W\), the norm of the tail after level \(N\) is at most \(m^{-(N+1)/2}\) by orthogonality. The remaining partial sum has finite rank. The first assertion follows from Theorem \ref{thm:hilbert}. For the second, take one of the \(h_Q\) from Proposition \ref{prop:noniso} at each level, and normalize it in \(L^1\). If \(Q\) has generation \(n-1\), then
 \[
   \norm{h_Q}_1=\mu(Q)^{1/2}\,\frac{1}{m}\sum_i\abs{w_i}.
 \]
 The images of the normalized vectors have a fixed positive \(L^2\) norm and are mutually orthogonal.
\end{proof}

\section{A rank-one criterion for the Dunford--Pettis property}
\label{sec:rankone}

Write
\[
  \Rank(W)=\Span\{v\otimes a\in W : v\in V_m,\ a\in\R^k\}.
\]
We first record the scalar product construction used in the positive direction and in the Hilbert-space argument for the example in Section \ref{sec:jordan}.

\begin{lemma}[Scalar Riesz products]
 \label{lem:riesz}
 Let \(U\subset V_m\) be a scalar constraint and let \(z_j\in L^1_U\) be blocks on successive disjoint intervals of levels with \(\sup_j\norm{z_j}_\infty\le M\), where \(M>0\). Let \(\sigma_j\) denote the coordinate signs on the independent product space \(\Sigma=\{-1,1\}^{\mathbb N}\). There is a contraction \(B:L^1(\Sigma)\to L^1_U\) with \(B\sigma_j=z_j/M\). If \(A:L^1_U\to H\) is bounded and \(H\) is a Hilbert space, then
 \[
   \sum_j\norm{Az_j}_H^2\le M^2\norm A^2.
 \]
\end{lemma}

\begin{proof}
 For each \(N\) define \(P_N(\sigma,t)=\prod_{j=1}^N(1+\sigma_j z_j(t)/M)\). For fixed signs this product is nonnegative and has integral one in \(t\), by conditioning successively at the start of each block. At a level in block \(j\), the later factors have conditional mean one after successive conditioning. The earlier factors are measurable at the start of block \(j\). Consequently the difference of \(P_N\) there is a predictable scalar multiple of \(d_n z_j\), so \(P_N(\sigma,\cdot)\in L^1_U\).

 For \(\phi\) depending on the first \(N\) signs, set \(B\phi=\E_\sigma[\phi P_N]\). The averages are consistent in \(N\). Since \(\norm{P_N(\sigma,\cdot)}_1=1\), \(B\) extends to a contraction. Averaging against \(\sigma_j\) gives \(B\sigma_j=z_j/M\). If \(H\) is Hilbert, apply Bessel's inequality in \(L^2(\Sigma;H)\) to \(\sigma\mapsto A(P_N(\sigma,\cdot))\) and the first \(N\) signs to obtain the asserted bound, then let \(N\to\infty\).
\end{proof}

We use the following consequence of the finite-rank conditional expectations and Lemma \ref{lem:stop}. If a bounded operator \(T:L^1_W\to Z\) fails complete continuity, there are uniformly bounded blocks \(u_j\) on successive disjoint intervals for which \(\inf_j\norm{Tu_j}>0\), where boundedness of the blocks is in \(L^\infty\). Indeed, start with a weakly null sequence \((f_j)\) for which the images are bounded away from zero. Inductively, given an endpoint \(a\), pass to an index with \(\norm{\E_a f_j}_1\) small, using finite rank. Then choose \(b>a\) with \(\norm{\E_b f_j-f_j}_1\) small. The blocks \((\E_b-\E_a)f_j\) can be made to differ from the selected \(f_j\) by errors tending to zero. They form a uniformly integrable family. Apply Lemma \ref{lem:stop} with a common error small enough that the images remain bounded away from zero.

\begin{theorem}
 \label{thm:criterion}
 The space \(L^1_W\) has the Dunford--Pettis property if and only if \(\Rank(W)=W\).
\end{theorem}

\begin{proof}
 First suppose that \(\Rank(W)=W\). The case \(W=0\) is finite-dimensional. In the scalar case, given a weakly compact \(T\), Lemma \ref{lem:riesz} gives \(Tu_j=M T B\sigma_j\) for any sequence of uniformly bounded blocks \(u_j\). The signs are weakly null in \(L^1(\Sigma)\), for example by orthonormality in \(L^2(\Sigma)\) and continuous inclusion into \(L^1(\Sigma)\). The operator \(TB\) is weakly compact, so the DPP of \(L^1(\Sigma)\) gives \(\norm{Tu_j}\to0\). The block reduction proves the scalar case.

 In the vector case choose a basis \(v_q\otimes a_q\), \(1\le q\le r\), of \(W\). For a sequence \((u_j)\) of blocks bounded in \(L^\infty\), expand every local difference in this basis and let \(u_{j,q}\) be the scalar block formed by the \(q\)-th coefficient times \(v_q\) at each parent. Then
 \[
   u_j=\sum_{q=1}^r a_q u_{j,q},\qquad
   u_{j,q}\in L^1_{\Span\{v_q\}},\qquad
   \norm{u_{j,q}}_2\le C_q\norm{u_j}_2 .
 \]
 The last inequality follows from boundedness of the local coordinate maps and orthogonality of differences in \(L^2\). For each \(q\), the scalar blocks are \(L^2\)-bounded and supported on disjoint intervals of levels, so they are weakly null in \(L^2\) and hence in \(L^1_{\Span\{v_q\}}\). For weakly compact \(T:L^1_W\to Z\), the map \(h\mapsto T(a_q h)\) is weakly compact on \(L^1_{\Span\{v_q\}}\). The scalar case shows that \(\norm{T(a_q u_{j,q})}\to0\). Summing over \(q\) and using the block reduction proves the DPP.

 Now suppose \(\Rank(W)\ne W\) and choose \(0\ne\phi\in W^*\) vanishing on \(\Rank(W)\). The local increase in the norm used below is the convexity defect from \cite[Sections~2--3]{ASW}. We use its telescoping sum to bound, uniformly on the unit ball, the number of coordinates of the detecting operator above any fixed threshold. For \(x\in\R^k\), \(w=(w_i)_{i=1}^m\in W\), set
 \[
    S(x,w)=\frac1m\sum_i\abs{x+w_i},
    \qquad D(x,w)=S(x,w)-\abs{x}.
 \]
 The triangle inequality gives \(D\ge0\). At equality, the child values \(x+w_i\) lie on a common ray or are all zero. Thus \(w\) has rank at most one and \(\phi(w)=0\). On the set \(S=1\), both \(x\) and the \(w_i\) are bounded, since \(x\) is the average of the child values. This set is compact. By continuity on that set and homogeneity, for every \(\epsilon>0\) there is a \(C_\epsilon\) with
 \begin{equation}
 \label{eq:defect}
  \abs{\phi(w)}\le\epsilon S(x,w)+C_\epsilon D(x,w).
 \end{equation}
 Indeed, on \(S=1\) the set \(\{\abs{\phi(w)}\ge\epsilon\}\), when nonempty, has \(D\) bounded away from zero.

 For a function \(f\), let \(w_n(t)=D_Qf\) on each parent \(Q\) of generation \(n-1\), and let \(D_n(t)=D(\E_{n-1}f(t),w_n(t))\). Define
 \[
   (Tf)_n=\E\phi(w_n).
 \]
 We have \(\E D_n=\norm{\E_nf}_1-\norm{\E_{n-1}f}_1\) and \(\E S=\norm{\E_nf}_1\le\norm f_1\) at that level. Thus \(\sum_n \E D_n\le\norm f_1\). By \eqref{eq:defect}, \(T\) is bounded into \(\ell^\infty\). For \(\norm f_1\le1\), at most \(C_\epsilon/\epsilon\) of its coordinates have absolute value exceeding \(2\epsilon\). This puts \(Tf\) in \(c_0\). The same bound on the number of large coordinates holds throughout the coordinatewise closure of \(T(B_{L^1_W})\), since exceeding a strict threshold in any fixed finite set of coordinates persists for sufficiently close points. The closure is therefore contained in \(c_0\). On a uniformly bounded subset of \(c_0\), weak and coordinatewise convergence coincide, by approximation of \(\ell^1=(c_0)^*\) functionals by finite sums. Compactness of the coordinatewise closure in the product of bounded intervals proves that \(T:L^1_W\to c_0\) is weakly compact.

 Choose \(w\in W\) with \(\phi(w)=1\) and let \(f_j(t)=w_{t_j}\), where \(t_j\) denotes the \(j\)-th digit of \(t\in\Omega\). The \(f_j\) form a bounded orthogonal sequence in \(L^2_W\) and so are weakly null in \(L^1_W\). The image \(Tf_j\) is the \(j\)-th unit vector of \(c_0\). This proves failure of the DPP.
\end{proof}

\section{Examples on the 4-ary tree}
\label{sec:examples}

\subsection{Three elementary constraints}

For every \(W\ne0\), repeating one pattern \(b=(b_i)\in W\) on independent levels gives an isomorphic copy of \(\ell^2\) in \(L^1_W\). Normalize \(b\) by \(m^{-1}\sum_i\abs{b_i}^2=1\) and set \(r_n(t)=b_{t_n}\). For a finitely supported scalar sequence \(c\), the sum \(u=\sum_n c_n r_n\) has \(\norm u_2=\norm c_2\) and \(\norm u_4\lesssim_b\norm c_2\), by expanding fourth moments of the independent centered coordinates. Interpolation gives \(\norm u_1\ge\norm u_2^3/\norm u_4^2\gtrsim_b\norm c_2\). The issue in the examples below is bounded complementation of a Hilbert copy.

Let \(v=(1,-1,0,0)\), \(w=(0,0,1,-1)\), and \(h=(1,1,-1,-1)\). The scalar space \(W_{\mathrm{s}}=\Span\{h\}\) satisfies (S2) and fails (S1). Theorem \ref{thm:criterion} gives the DPP, so there is no complemented \(\ell^2\). As in \cite[Section~6.1]{KL}, when each \(4\)-adic digit is read as two binary digits, \(h\) corresponds to the Haar difference that is revealed by the first digit in the pair. Thus \(L^1_{W_{\mathrm{s}}}\) is the constants plus the closed span of Haar functions on dyadic intervals of lengths \(2^{-2n}\), \(n\ge0\). Speckhofer \cite{Speckhofer} proved nonisomorphism with \(L^1\) for the closed span of these levels.

For the Hilbertian example, take the space \(W_{\mathrm{H}}=\Span\{B\}\) of \(4\times2\) matrices, with
\[
 B=\begin{pmatrix}
 1&0\\
 0&1\\
 -1&-1\\
 0&0
 \end{pmatrix}.
\]
It satisfies (S1). In fact a complemented \(\ell^2\) can be exhibited for any nonzero \(W\) satisfying (S1). Choose \(b=(b_i)_{i=1}^m\in W\) with \(m^{-1}\sum_i\abs{b_i}^2=1\), and set \(r_n(t)=b_{t_n}\). The maps
\[
   J:\ell^2\to L^1_W,\quad Jc=\sum_n c_n r_n ,
   \qquad
   Af=(\E\langle f,r_n\rangle)_{n\ge1}
\]
are bounded. For \(J\) this follows from orthonormality in \(L^2_W\). For \(A\), \(\E\langle f,r_n\rangle=\E\langle d_nf,r_n\rangle\), since \(r_n\) is \(\F_n\)-measurable with conditional mean zero at time \(n-1\). The triangle inequality in \(\ell^2\), followed by \eqref{eq:Janson}, gives
\[
  \norm{Af}_2\le \max_i\abs{b_i}\,
       \E\left(\sum_n\abs{d_nf}^{2}\right)^{1/2}
       \lesssim_W\norm f_1 .
\]
Orthonormality yields \(AJ=I_{\ell^2}\). Thus \(JA\) is a projection onto \(J(\ell^2)\), and \(L^1_W\) fails the DPP.

One can also obtain a complemented \(\ell^2\) under (S2) with (S1) failing. Take
\[
  W_{\mathrm{split}}
   =\{\, (aB,bh): a,b\in\R\,\}\subset V_4\otimes\R^3,
\]
where \(bh\) is the third column. If \(a\ne0\), the matrix has rank at least two. For \(a=0\), the only nonzero rank-one pattern is a multiple of \(h\) in the third coordinate, which is allowed by (S2). The projections on the first two coordinates and on the third coordinate yield \(L^1_{W_{\mathrm{split}}}\simeq L^1_{W_{\mathrm H}}\oplus L^1_{W_{\mathrm s}}\). The first summand supplies a complemented \(\ell^2\). This coordinate splitting is related to the \(L^1\oplus H^1\) example of Kazaniecki and Lechner \cite[Theorem~6.6]{KL}. Here both summands satisfy (S2).

\subsection{A Jordan-type constraint}
\label{sec:jordan}

Define
\[
 W_{\mathrm J}
 = \{\,(av+bw,bv):a,b\in\R\,\}\subset V_4\otimes\R^2 .
\]
The name comes from the shear \(I+N\) on the value space, where \(N(x,y)=(y,0)\). Acting on each child value, \(N\) sends the generator \((w,v)\) to \((v,0)\) and annihilates \((v,0)\). Thus \(I+N\) acts on \(W_{\mathrm J}\) as a single Jordan block. The same map supplies the admissible correction \((V_j,0)\) in the proof below.

The four row vectors are \((a,b),(-a,-b),(b,0),(-b,0)\). The determinant of the first and third rows is \(-b^2\), so the rank-one span is given by \(b=0\). The only nonzero rank-one patterns have first column a multiple of \(v\) and second column zero. In particular \(W_{\mathrm J}\) satisfies (S2) and fails (S1).

\begin{theorem}
 \label{thm:jordan}
 The space \(L^1_{W_{\mathrm J}}\) fails the Dunford--Pettis property and contains no infinite-dimensional complemented reflexive subspace. Every bounded operator from \(L^1_{W_{\mathrm J}}\) to a Hilbert space is completely continuous.
\end{theorem}

The criterion in Theorem \ref{thm:criterion} gives failure of the DPP. The following argument gives an explicit witness with a reflexive target, \(\ell^4\).

\begin{proof}[A witness into \(\ell^4\)]
 Let \(x\) be a parent value and \(d_i\) its four increment rows. As above put \(S=\tfrac14\sum_i\abs{x+d_i}\) and \(D=S-\abs x\). For a unit vector \(e\) in the direction of \(x\), choosing any \(e\) when \(x=0\), we have
 \[
   \sum_i(\abs{x+d_i}-\langle e,x+d_i\rangle)=4D .
 \]
 Since \(\abs{x+d_i}\le4S\), the component of \(d_i\) perpendicular to \(e\) has modulus at most \((32SD)^{1/2}\). This follows by factoring
 \(\abs{x+d_i}^2-\langle e,x+d_i\rangle^2\). Also \(\abs{d_i}\le5S\). Estimating the determinant of \(d_1,d_3\) in these coordinates gives \(\abs b^2\lesssim S^{3/2}D^{1/2}\).

 On level \(n\) write \(b_n\) for \(b\) as a function of the parent. By H\"older's inequality and the telescoping estimate for the defects,
 \[
   \abs{\E b_n}^4\lesssim \norm f_1^3\,\E D_n,\qquad
   \sum_{n\ge1}\abs{\E b_n}^4\lesssim \norm f_1^4 .
 \]
 Therefore \(f\mapsto(\E b_n)_n\) is a bounded operator into \(\ell^4\). The functions \(f_j(t)=(w_{t_j},v_{t_j})\) are weakly null by \(L^2\)-orthogonality and their images are unit vectors. The operator is weakly compact because \(\ell^4\) is reflexive.
\end{proof}

\begin{lemma}[Survival estimate]
 \label{lem:survival}
 Suppose \(Z_0=1\) and \(Z_0,\ldots,Z_n\) is a real martingale with increments of modulus at most \(1/2\), absorbed when it first reaches \((-\infty,0]\). Then
 \[
   \mathbb P(Z_n>0)\ge\frac{1}{4\sqrt{n+1}} .
 \]
\end{lemma}

\begin{proof}
 Put \(a=1+\sqrt{n+1}\) and stop on leaving \((0,a)\) or at time \(n\), whichever comes first. The stopped value has expectation \(1\). It is nonpositive on the event of exit below and at most \(a+1/2\) otherwise. Thus the probability of exiting above or staying inside to time \(n\) is at least \((a+1/2)^{-1}\). Conditional on a first exit above at time \(\tau\), the chance of a subsequent fall to zero before time \(n\) is at most \(n/(4a^2)\le1/4\), by applying the maximal inequality to the square of the martingale of increments after \(\tau\). The conditional second moment of the sum of those increments is at most \(n/4\). Consequently
 \[
   \mathbb P(Z_n>0)\ge \frac{3}{4(a+1/2)}
                \ge\frac{1}{4\sqrt{n+1}} .
 \]
\end{proof}

\begin{proof}[Proof of the Hilbert-space assertion of Theorem \ref{thm:jordan}]
 Suppose \(A:L^1_{W_{\mathrm J}}\to H\), with \(H\) Hilbert, fails complete continuity. Choose uniformly \(L^\infty\)-bounded successive blocks \(U_j=(K_j,V_j)\) as in the block reduction, with \(\inf_j\norm{AU_j}>0\). Let \(M\ge1\) bound \(\abs{K_j}\) and \(\abs{V_j}\). By the Jordan constraint, \(V_j\) is a scalar block in \(L^1_{\Span\{v\}}\) on the same interval. Hence \((V_j,0)\) is an allowed block. Fix \(0<\eta\le 1/(4M)\).

 The construction gives the \(j\)-th block a sign coefficient of order \(j^{-1/2}\), up to a scalar correction that is square summable after applying \(A\). The lower bound of order \(j^{-1/2}\) suffices because its squares have divergent sum.

 We construct, for each sign sequence \(\sigma\), nonnegative coordinates \((F_j,G_j)\) in \(L^1_{W_{\mathrm J}}\) with initial constant value \((F_0,G_0)=(2,1)\). The update at block \(j\) takes place on the event
 \[
     C_j=\{F_{j-1}/G_{j-1}>1\}.
 \]
 We will have \(G_j>0\). With \(\theta_j\in[0,1]\) defined below, put
 \begin{equation}
 \label{eq:Jordanupdate}
 \begin{pmatrix}F_j\\ G_j\end{pmatrix}
 =
 \begin{pmatrix}F_{j-1}\\ G_{j-1}\end{pmatrix}
 +\eta\sigma_j\theta_j\mathbf 1_{C_j}
 \left[
    G_{j-1}\begin{pmatrix}K_j\\ V_j\end{pmatrix}
   +F_{j-1}\begin{pmatrix}V_j\\0\end{pmatrix}
 \right].
 \end{equation}
 Here \(\theta_j\) depends on the tree variable only and is measurable at the start of block \(j\). Both coordinates of the previous pair and \(\mathbf 1_{C_j}\) are also measurable there. Thus \eqref{eq:Jordanupdate} adds an allowed block for each \(\sigma\). The two coordinates stay nonnegative. Indeed, \(G_j\) is obtained by multiplying \(G_{j-1}\) by a positive factor. On \(C_j\), the correction in the first coordinate cancels the term proportional to \(R_{j-1}V_j\) when we form the ratio \(R_j=F_j/G_j\). In fact,
\[
 R_j-R_{j-1}
   =\frac{\eta\sigma_j\theta_j K_j}{1+\eta\sigma_j\theta_j V_j}.
\]
Its modulus is at most \(\eta M/(1-\eta M)\le 1/2\). When \(C_j\) fails, both coordinates are frozen. Since the update has zero mean in the tree variable, it follows that
 \begin{equation}
 \label{eq:Jordanbound}
    \norm{(F_j,G_j)}_1\le 3 \qquad\text{for each sign sequence}.
 \end{equation}

 To choose the predictable multiplier, hold the tree variable \(t\) fixed. Both \(F_i\) and \(G_i\) are martingales in the first \(i\) signs, by the factor \(\sigma_i\) in their increments. Through time \(j-1\), use the sign probability with density \(G_{j-1}\) relative to fair signs. It has total mass one and its density through time \(i\le j-1\) is \(G_i\). Since \(G_i R_i=F_i\), the ratio \(R_i\) is a martingale under this probability. The process \(R_i-1\) starts at one, has increments of modulus at most \(1/2\), and is absorbed on reaching \((-\infty,0]\). By Lemma \ref{lem:survival}, the probability of being active at step \(j\) satisfies
 \[
   H_j(t):=\E_\sigma[G_{j-1}(t)\mathbf 1_{C_j}(t)]
               \ge \gamma_j,\qquad \gamma_j=\frac{1}{4\sqrt j}.
 \]
 Set \(\theta_j=\gamma_j/H_j\). The bound gives \(0<\theta_j\le1\). Since \(H_j\) is determined before block \(j\), this makes \eqref{eq:Jordanupdate} a recursion. Its other averaged coefficient satisfies
 \[
  0\le \beta_j:=\theta_j\E_\sigma[F_{j-1}\mathbf 1_{C_j}]\le 2,
 \]
 using nonnegativity and the sign-martingale mean of \(F_{j-1}\).

 For \(j\le N\), the coefficient of \(\sigma_j\) in \((F_N,G_N)\), computed by averaging over signs, is
 \[
   \E_\sigma[\sigma_j(F_N,G_N)]
     =\eta\left(\gamma_j U_j+(\beta_j V_j,0)\right).
 \]
 Only the update at time \(j\) enters this coefficient, by the sign-martingale property. Applying \(A\), Bessel's inequality in the signs and \eqref{eq:Jordanbound} give
 \[
   \sum_{j=1}^N
   \norm{\gamma_j AU_j+A(\beta_j V_j,0)}_H^2
   \le 9\eta^{-2}\norm A^2 .
 \]
 The terms \(A(\beta_j V_j,0)\) are square summable, by Lemma \ref{lem:riesz} applied to the uniformly bounded scalar blocks \(\beta_j V_j\) and the operator \(z\mapsto A(z,0)\). Therefore \(\sum_j\gamma_j^2\norm{AU_j}_H^2<\infty\), which contradicts the lower bound on \(\norm{AU_j}_H\) and the divergence of \(\sum_j \gamma_j^2\). This proves complete continuity for operators into \(H\). Both complementability conclusions now follow from Theorem \ref{thm:hilbert}.
\end{proof}

\begin{table}[htbp]
\centering
\begin{tabular}{@{}lcccc@{}}
\toprule
Constraint on the \(4\)-ary tree & (S1) & (S2) & DPP & Complemented \(\ell^2\)\\
\midrule
\(W_{\mathrm s}=\Span\{h\}\) & No & Yes & Yes & No\\
\(W_{\mathrm H}=\Span\{B\}\) & Yes & Yes & No & Yes\\
\(W_{\mathrm{split}}=\{(aB,bh):a,b\in\R\}\) & No & Yes & No & Yes\\
\(W_{\mathrm J}=\{(av+bw,bv):a,b\in\R\}\) & No & Yes & No & No\\
\bottomrule
\end{tabular}
\caption{All three possible combinations occur under (S2), already among constraints that fail (S1).}
\label{tab:examples}
\end{table}

\section{Comparison with Hardy and Sobolev spaces}
\label{sec:comparison}

\subsection{Hardy spaces}

For analytic \(H^1(\mathbb D)\), the coefficient map
\[
 f\longmapsto \bigl(\widehat f(2^j)\bigr)_{j\ge1}
\]
is bounded into \(\ell^2\) by Paley's theorem. Synthesis by the corresponding monomials is bounded from \(\ell^2\) to \(H^1\) by the inclusion of the boundary \(L^2\) space in \(L^1\). The composition is the identity on \(\ell^2\). Thus analytic \(H^1\) fails the DPP, and the theorem of Kwapie\'n and Pe{\l}czy\'nski characterizes all its complemented reflexive subspaces as Hilbertian \cite{KP}. Maurey's isomorphism between analytic \(H^1\) and dyadic martingale \(H^1\) gives the corresponding statement for the dyadic space \cite{Maurey}. The coefficient projection under (S1) in Section \ref{sec:examples} is a direct martingale construction. Theorem \ref{thm:hilbert} extends the Hilbertian conclusion to constraints for which the Hardy norm estimate fails.

\subsection{Isotropic Sobolev spaces}

On a bounded interval \(I\), evaluation at a fixed point and differentiation identify \(W^{1,1}(I)\) with \(\R\oplus L^1(I)\). On \(\T\), the mean and derivative identify \(W^{1,1}(\T)\) with \(\R\oplus L^1_0(\T)\), with the subscript denoting zero mean. Both spaces have the DPP by the classical \(L^1\) result.

In dimensions \(d\ge2\), the failure of the DPP for the homogeneous spaces \(\dot W^{1,1}(\Omega)\) and \(\dot{BV}(\Omega)\) was announced by Nazarov, reporting joint work with Kazaniecki and Wojciechowski, at an Ohio State University seminar on 24 September 2026 \cite{KNW}. We refer to the public seminar abstract for the announcement. The earlier paper of Kazaniecki, Tselishchev and Wojciechowski \cite{KTW} relates the DPP for \(W^{1,1}\) and \(SBV\) through the Schur property of jump parts of gradients. On a cube, Poincar\'e's inequality identifies the homogeneous Sobolev space modulo constants with the mean-zero subspace of the inhomogeneous space. A failure of the DPP on a cube transfers to \(W^{1,1}(\T^d)\) and \(W^{1,1}(\R^d)\) by complementation: restriction to the cube has bounded linear right inverses obtained by reflection, with a cutoff in the Euclidean case.

\begin{theorem}
\label{thm:torus}
For every integer \(n\ge1\), each complemented reflexive subspace of
\(W^{1,1}(\T^n)\) is isomorphic to a Hilbert space.
\end{theorem}

\begin{proof}
We verify \eqref{eq:approx} with \(X=W^{1,1}(\T^n)\), \(H=W^{1,2}(\T^n)\), and the natural inclusion. The proof uses Lipschitz truncation in place of martingale stopping.
Let \(K\) be weakly compact in \(W^{1,1}(\T^n)\), write \(M_{\rm HL}\) for the periodic Hardy--Littlewood maximal operator, and set
\[
  E_\lambda(u)=\{M_{\rm HL}(\abs{\nabla u})>\lambda\},
                       \qquad u\in K .
\]
The family of gradients is uniformly integrable. The weak \((1,1)\) estimate for the maximal operator applied to \(\abs{\nabla u}\mathbf1_{\{\abs{\nabla u}>\lambda/2\}}\) yields
\[
  \lambda|E_\lambda(u)|
     \lesssim_n \int_{\{\abs{\nabla u}>\lambda/2\}}\abs{\nabla u}
       \longrightarrow 0
\]
uniformly over \(K\). On the complement of \(E_\lambda(u)\) the Sobolev representative is \(C_n\lambda\)-Lipschitz outside a null set, by the standard pointwise estimate in terms of \(M_{\rm HL}(\abs{\nabla u})\) at the two points. Extend it to a function \(u_\lambda\) on the torus with the same Lipschitz bound, as in the classical Lipschitz truncation argument \cite{AF,DMS}. Locality of the weak gradient then gives
\[
   \norm{\nabla u-\nabla u_\lambda}_1
    \le \int_{E_\lambda(u)}\abs{\nabla u}+C_n\lambda|E_\lambda(u)| .
\]
The right side tends to zero uniformly in \(K\). Adjust the mean of \(u_\lambda\) to that of \(u\). Poincar\'e's inequality gives uniform \(W^{1,1}\) approximation from the gradient estimate, and gives a uniform \(W^{1,2}\) bound for the approximants from the Lipschitz bound and the bounded means. This yields \eqref{eq:approx}, and Lemma~\ref{lem:hilbert-approx} applies.
\end{proof}

The conclusion also holds on \(\R^n\). Indeed, a weakly compact set in
\(W^{1,1}(\R^n)\) has uniformly small \(L^1\) tails for both the functions
and their gradients, by the weak compactness criterion in \(L^1\).
Given an approximation error, first multiply by a common smooth cutoff
equal to one on a sufficiently large ball. Regard the cut-off functions
as periodic on a cube containing their common support in its interior.
Apply the preceding approximation there, then multiply the approximants
by a second smooth cutoff equal to one on that support and extend by zero.
This gives uniformly bounded \(W^{1,2}\) approximants with common compact
support. For the unit ball of a complemented reflexive subspace, choose the
error so that the projected error is at most \(1/2\). The same Hilbert
space of functions with fixed support then suffices at each step of the
geometric iteration in Lemma~\ref{lem:hilbert-approx}.

Consequently the conclusion transfers to every open set \(\Omega\subset\R^n\)
for which restriction from \(W^{1,1}(\R^n)\) has a bounded linear right
inverse, in particular to bounded Lipschitz domains. In the plane there is
a further transfer without boundary regularity. In joint work with
Derezi\'{n}ski and Nazarov, Wojciechowski announced that for every bounded
simply connected domain \(\Omega\subset\R^2\) and \(1\le p<2\), the
spaces \(W^{1,p}(\Omega)\) and \(W^{1,p}((0,1)^2)\) are isomorphic, as
are their homogeneous counterparts \cite{DNW}. At \(p=1\), this announced classification
transfers the Hilbertian conclusion from the square to all such domains.

Kazaniecki, Pakosz and Wojciechowski \cite[Theorem~1]{KPW} prove that an operator on \(L^2(\T^d)\) factoring through \(W^{1,1}(\T^d)\) by two translation-invariant maps belongs to the Schatten class \(\mathcal S_p\) for every \(p>2d+4\). This excludes an infinite-rank orthogonal projection obtained by such a factorization.

The critical embedding \(W^{1,1}(\T^2)\hookrightarrow L^2(\T^2)\) is completely continuous by \eqref{eq:approx} and the compact embedding \(W^{1,2}(\T^2)\hookrightarrow L^2(\T^2)\). Its failure of absolute \(2\)-summability is Kislyakov's result \cite{Kislyakov}. Recent factorization and summability estimates for this embedding are given in \cite{KPaW}. These properties parallel those of \(I_{1/2}\) under (S2), including in the Jordan example.

\subsection{Anisotropic Sobolev spaces}

For a finite set \(S\subset\mathbb N_0^d\) containing zero, let \(L^1_S(\T^d)\) denote the completion of trigonometric polynomials with respect to the norm \(\sum_{\alpha\in S}\norm{\partial^\alpha u}_1\). Pe{\l}czy\'nski and Wojciechowski \cite{PW} characterize the existence of a translation-invariant complemented Hilbert subspace of infinite dimension in these spaces. Their criterion is that two points of \(S\) with opposite parity of derivative order lie on a supporting hyperplane of \(\operatorname{conv}(S)\) which is parallel to no coordinate axis, or that this condition holds after restricting \(S\) to a coordinate subspace and working in that subspace. For \(S=\{(0,0),(1,0),(0,2)\}\) the line through \((1,0)\) and \((0,2)\) gives such a hyperplane. The resulting space, with norm
\[
   \norm u_1+\norm{\partial_x u}_1+\norm{\partial_y^2 u}_1 ,
\]
therefore contains a complemented copy of \(\ell^2\) and fails the DPP.
For the isotropic first-order set \(S=\{0,e_1,\ldots,e_d\}\), a supporting hyperplane through zero and any \(e_i\) is parallel to the \(i\)-th axis. The coordinate restrictions have the same obstruction. The criterion consequently excludes translation-invariant complemented Hilbert subspaces of infinite dimension in this case. See \cite{WojProj} for further results on translation-invariant projections. The Paley projection criterion proves failure of the DPP in the cases where it supplies a complemented Hilbert space. In the other cases it leaves the DPP undecided, as the isotropic example illustrates.

\section*{Use of Artificial Intelligence}

The results in this article were developed in collaboration with
several different AI models. The author has carefully read the entire
final manuscript and takes full responsibility for its contents.

\end{document}